\documentclass[sn-mathphys-num]{sn-jnl}
\usepackage{geometry} 
\usepackage{graphicx}%
\usepackage{multirow}%
\usepackage{amsmath,amssymb,amsfonts,placeins,stfloats}%
\usepackage{amsthm}%
\usepackage{mathrsfs}%
\usepackage[title]{appendix}%
\usepackage{xcolor}%
\usepackage{textcomp}%
\usepackage{manyfoot}%
\usepackage{booktabs}%
\usepackage{algorithm}%
\usepackage{algorithmicx}%
\usepackage{algpseudocode}%
\usepackage{listings}%
\usepackage{epsfig} 
\usepackage{epstopdf}
\usepackage{hyperref}
\usepackage{sidecap}
\usepackage{natbib}
\usepackage{wrapfig}
\usepackage{subcaption}
\usepackage{float}
\usepackage{lmodern}

\makeatletter\def\@font@warning#1{}\makeatother

\theoremstyle{thmstyleone}%
\newtheorem{theorem}{Theorem}
\theoremstyle{thmstyletwo}%

\theoremstyle{thmstylethree}%
\newtheorem{definition}{Definition}%
\newtheorem{lem}{Lemma}

\begin{document}

\title[Functional $\mathcal{H}_{\infty}$ Filtering for Descriptor Systems with Monotone nonlinearities]{Functional $\mathcal{H}_{\infty}$ Filtering for Descriptor Systems with Monotone nonlinearities}


\author[1]{\fnm{Rishabh} \sur{Sharma}}\email{rishabh$\_$2021ma18@iitp.ac.in}
\equalcont{These authors contributed equally to this work.}

\author[2]{\fnm{Mahendra Kumar} \sur{Gupta}}\email{mkgupta@iitbbs.ac.in}
\equalcont{These authors contributed equally to this work.}

\author*[1]{\fnm{Nutan Kumar} \sur{Tomar}}\email{nktomar@iitp.ac.in}
\equalcont{These authors contributed equally to this work.}

\affil[1]{\orgdiv{Department of Mathematics}, \orgname{Indian Institute of Technology Patna}, \orgaddress{ \city{Patna}, \state{Bihar}, \postcode{801106}, \country{India}}}

\affil[2]{\orgdiv{School of Basic Sciences}, \orgname{Indian Institute of Technology Bhubaneswar}, \orgaddress{ \city{Bhubaneswar}, \state{Odisha}, \postcode{752050}, \country{India}}}


\abstract{This paper introduces a novel approach to design of functional $\mathcal{H}_{\infty}$ filters for a class of nonlinear descriptor systems subjected to disturbances. Departing from conventional assumptions regarding system regularity, we adopt a more inclusive approach by considering general descriptor systems that satisfy a rank condition on their coefficient matrices. Under this rank condition, we establish a linear matrix inequality (LMI) as a sufficient criterion ensuring the stability of the error system and constraining the $\mathcal{L}_2$ gain of the mapping from disturbances to errors to a predetermined level. The efficacy of the proposed approach is demonstrated through a practical example involving a simple constrained mechanical system.}




\maketitle

\section{Introduction}\label{intro}
Mathematical models play a crucial role in understanding and analyzing control systems. Control theory often deals with physical systems represented by ordinary differential equations (ODEs) along with certain algebraic constraints. Such systems are called descriptor systems, which are also known as generalized state-space systems, implicit systems, singular systems, or systems described by differential algebraic equations (DAEs). In some cases, it is feasible to solve the algebraic constraints explicitly, transforming the system model into a well-known state space model, which is essentially a set of ODEs. Nevertheless, these transformations necessitate human intervention, adjustments or eliminations of system variables. As a result, the transformed state space models do not capture many interesting intrinsic properties of the underlying physical phenomena. For example, state space models can not effectively handle impulses occurring in electrical circuits \cite{belur2019persistence}. Consequently, it becomes practically important to study the properties of descriptor systems. For a more in-depth exploration of descriptor systems and their applications in engineering, readers are encouraged to see the books \cite{eich1998numerical, riaza2008differential, kumar1999control} and the literature cited therein. 

The focus of this paper is on state estimation for descriptor systems. State estimation is a classical problem in control systems theory, and the first solution to this problem for linear time-invariant (LTI) state space systems dates back to the early 1960s by Luenberger \cite{luenberger1964observing} and Kalman \cite{kalman1960new}. The Luenberger approach gave rise to the well-known observer design problem. A great deal of work has been devoted to the observer design problem for linear descriptor systems, and for a thorough discussion on the existence conditions, we refer the readers to Jaiswal \emph{et al.}    \cite{jaiswal2021necessary}. Moreover, in the context of linear descriptor systems with corrupted input and output data, a significant amount of research has focused on state estimation using Kalman filtering approach  \cite{dai1989singular,ishihara2005optimal, goel2024kalman}. However, it is notable that the Kalman filtering methods assume that the plant system has dynamics described by external noises with precisely known statistical characteristics \cite{brown1997introduction}. To overcome these limitations, an alternative approach, $\mathcal{H}_{\infty}$ filtering, was first proposed by Xu \emph{et al.} \cite{xu2003h} for square and regular linear descriptor systems. This method does not rely on assuming any statistical characteristics for the exogenous noises, but only requires bounded (within the $\mathcal{L}_2$-norm) noises   \cite{xu2007reduced,gao2020robust, tunga2023h}. Moreover, $\mathcal{H}_{\infty}$ filtering proves to be more robust than Kalman filtering, especially when dealing with additional uncertainties in descriptor systems \cite{xu2006robust}.

In the present paper, we study the $\mathcal{H}_{\infty}$ filtering problem for a class of nonlinear descriptor systems. Nonlinear descriptor systems seem to have been first considered by Luenberger \cite{luenberger1979non}. Thereafter, several results on nonlinear descriptor systems are available in the literature, and for solutions to such systems, we refer the readers to the articles \cite{petzold1986numerical,liberzon2012switched} and the books \cite{kunkel2006differential,lamour2013differential}. It is important that the research on state estimation for nonlinear descriptor systems emphasizes on systems characterized by structured nonlinearities. The observer design problem for descriptor systems with Lipschitz nonlinear functions has been the topic of considerable research over the past three decades, see \cite{ shields1997observer,  lu2006full,gupta2014observer,berger2018observers} and the references therein. Roughly speaking, a Lipschitz function is characterized by the condition that the absolute value of its slope is bounded above by a Lipschitz constant. On the other hand, monotone nonlinearities encompass functions that exhibit lower bounds on their slope. Notably, the monotone nonlinearities arise in the modeling of many physical systems like stiffening springs in constrained mechanical systems and capacitors in circuit systems \cite{fan2003observer,campbell1982singular}. Yang \emph{et al.} \cite{yang2013nonlinear} have introduced a full-order observer design for descriptor systems that meet a broader monotone condition, but the design is limited to square systems only. Meanwhile, Gupta \emph{et al.} \cite{gupta2018unknown} have proposed a reduced-order observer design tailored for nonlinear descriptor systems with generalized monotone nonlinearities.  Subsequently, in 2020, Berger \emph{et al.} \cite{berger2020observers} have devised an observer design framework for descriptor systems characterized by Lipschitz or monotone nonlinearities. Nonetheless, the observers in \cite{yang2013nonlinear, berger2020observers} adopt a descriptor form, which is not favored due to its implicit nature and the need for consistent initial conditions during simulation. Moysis \emph{et al.} \cite{moysis2020observer} expanded upon the work in \cite{gupta2018unknown}, focusing on nonlinearities that adhere to incremental quadratic constraints in output equations. The $\mathcal{H}_{\infty}$ filters for descriptor systems have been initially proposed for systems having Lipschitz nonlinearity \cite{darouach2011h} and Monotone-type nonlinearity \cite{yang2012h,aliyu2012cal} in the early 2010s.  Recently, state and adaptive disturbance $\mathcal{H}_{\infty}$ observer for a class of nonlinear descriptor systems with the disturbance generated by an unknown exogenous system is addressed in \cite{liu2023state}.

The aforementioned papers primarily focus on estimators designed to estimate the entire state vector of a given control system. However, in numerous applications, it is often sufficient to have information about only a subset or combination of states rather than the entire state vector. An estimator that targets a specific combination of states without estimating the entire state vector is called a functional (or partial) estimator. Such estimators find utility in various applications, particularly those involving feedback control, fault detection, and disturbance estimation \cite{trinh2011functional}. With the use of functional observers, $\mathcal{H}_{\infty}$ controller design challenges have been explored for both linear descriptor systems \cite{xu2006robust} and nonlinear descriptor systems \cite{wang2006h}. A primary theoretical objective in the study of functional estimators is to determine the conditions under which they can exist and offer accurate estimates. Notably, functional estimators can be designed with more relaxed assumptions compared to full-state estimators; for recent advancements in functional observer design for linear descriptor systems, the interested readers are refereed to \cite{jaiswal2021existence,jaiswal2022functional,tunga2023functional,jaiswal2024existence} and the references therein. To the best of our knowledge, there is no literature on functional $\mathcal{H}_{\infty}$ filters or observers specifically tailored for monotone nonlinear descriptor systems.

The current paper introduces a novel concept of the functional $\mathcal{H}_{\infty}$ filter for nonlinear descriptor systems with generalized monotone nonlinearities. The basis of our definition for functional $\mathcal{H}_{\infty}$ filter (cf. Definition $2$) relies on the behavioral solution theory for descriptor systems. The filter is designed to furnish a stable estimation for a given functional vector of semistates in the presence of disturbances. A key innovation of our proposed design methodology lies in its ability to provide the filter with an order that equal to the dimension of the functional vector. Notably, our approach refrains from making any assumptions about the system, such as regularity or squareness, thereby enhancing its applicability. Furthermore, the stability conditions for the error system is expressed using the Lyapunov approach in the form of an LMI.

The paper is organized as follows: Section \ref{prob} provides the problem statement and lays out the necessary preliminaries for our analysis. In Section \ref{main}, we establish the primary findings of the paper. Initially, a method is developed for determining coefficient matrices for the proposed filter, followed by the development of error system stability using the Lyapunov theory. Section \ref{num} presents a numerical example that serves to illustrate the theoretical results obtained. Lastly, Section \ref{concl} concludes the paper. 

\subsection{Notations}
  $$\begin{array}{ll}
			\mathbb{N}, \mathbb{N}_0 &  \mathbb{N}\text { is the set of natural numbers, and\,} \mathbb{N}_0 := \mathbb{N} \cup \{0\}; \\
           0, I & \text{The zero and identity matrices of appropriate dimensions, respectively;}\\
          \mathbb{C}, \mathbb{\overline{C}}^+  & \mathbb{C} \text{ is the set of complex numbers, and \,} \mathbb{\overline{C}}^+ := \{\alpha \mid \alpha \in \mathbb{C}, \; Re(\alpha) \geq 0\};\\
           A \in \mathbb{R}^{m \times n} & \text{The matrix A is in the set of real}\; m \times n \; \text{matrices;}\\
			A>\mathscr{} 0 &  \forall x \in \mathscr{V} \backslash\{0\},  \text { we have } x^{T} A x>0, A \in \mathbb{R}^{n \times n} \text { and } \mathscr{V} \text { is a subspace in } \mathbb{R}^n ; \\
	    	 A^{T}, A^{+} & \text {Transpose and Moore-Penrose inverse (MP-inverse) of } A \in \mathbb{R}^{m \times n}, \text { respectively}; \\
			\mathscr{C}^k(X \rightarrow Y) & \text { Set of } k \text {-times continuously differentiable functions } f: X \rightarrow Y,\; k \in \mathbb{N}_0;\\

   \|x(t)\| & \sqrt{x^T(t)x(t)}, \text{the Euclidean norm of}\;  x(t) \in \mathbb{R}^n;\\
   \mathcal{L}_2\big([0,t_f], \mathbb{R}^n\big) & \{f: \int_0^{t_{f}} \|f(t)\|^2 \mathrm{~d} t<\infty \; \text{for a fixed} \;  t_f>0 \};  \\  
    \text{dom}(f) &  \text{ The domain of the function} \; f; \\
   \mathbb{R}[\mu] & \;  \text{The ring of polynomials with coefficients in } \mathbb{R}.
    	\end{array}
	$$
		
\section{Problem Description and Preliminaries}\label{prob}
\noindent Consider a nonlinear descriptor system
			\begin{subequations}\label{sys}
			\begin{eqnarray}
				E\dot{x}(t)&= & Ax(t)+Bu(t)+Dv(t)+Fg(HKx,u),\label{sysa}\\
				y(t)&= & Cx(t)+Gv(t), \label{sysb}\\
				z(t)&= & Kx(t), \label{sysc}
			\end{eqnarray}
		\end{subequations}
		   where $x(t) \in \mathbb{R}^{n}$ is the semistate vector, $u(t) \in \mathbb{R}^{k}$ is the (known) control  input vector, $y(t) \in \mathbb{R}^{r}$ is the (measured) output vector, $z(t) \in \mathbb{R}^{p}$ is the (unmeasured) output functional vector, and $v(t) \in \mathbb{R}^{q}$ is the unknown disturbance vector. The coefficients $E, A \in {\mathbb{R}}^{m \times n}$, $B \in {\mathbb{R}}^{m \times k}$, $C \in {\mathbb{R}}^{r \times n}$, $D \in {\mathbb{R}}^{m \times q}$, $F \in {\mathbb{R}}^{m \times l}$,  $G \in {\mathbb{R}}^{r \times q}$, $H \in {\mathbb{R}}^{l \times p}$ and $K \in {\mathbb{R}}^{p \times n}$ are known constant matrices. Moreover, for some open sets $\mathscr{X}_1 \subseteq \mathbb{R}^p$, $\mathscr{U} \subseteq \mathbb{R}^k$ and $\mathscr{X}_2 = H\mathscr{X}_1 \subseteq \mathbb{R}^l$, the nonlinear function $g \in \mathscr{C}^1(\mathscr{X}_2 \times \mathscr{U} \rightarrow \mathbb{R}^l)$ satisfies a generalized monotone condition \cite{berger2020observers}
  \begin{eqnarray}\label{ineq}
     \forall \; x_{1}, x_{2} \in  \mathscr{X}_2 : (x_{1}-x_{2})^T\big(g(x_{1},u)-g(x_{2},u)\big) \geq \frac{1}{2} \rho \|x_{1}-x_{2}\|^2,
   \,  \text{where} \; \rho \in \mathbb{R}.
  \end{eqnarray}
   Clearly, if we define $\Delta g  =  g(x_1,u)-g(x_2,u)$ and $\Delta x = x_{1}-x_{2}$, then \eqref{ineq} reduces to 
  \begin{eqnarray}\label{ineq2}
      \Delta x^{T}\Delta g + \Delta g^{T} \Delta x \geq \rho \|\Delta x\|^2 . 
  \end{eqnarray}
  The property \eqref{ineq} (equivalently \eqref{ineq2}) is also called multivariable sector property, and it can be proved easily that \eqref{ineq} holds if the function $g$ satisfies the slope bound property (Lemma $1$ in \cite{liu2012unknown}):
   \begin{eqnarray} \label{ineq3}
		\forall \; s \in \mathbb{R}^{l}: \frac{\partial g(s)}{\partial s}+\left(\frac{\partial g(s)}{\partial s}\right)^T  \geq  \rho I_l.
		\end{eqnarray}
  Notably, if $\rho=0$, characterization \eqref{ineq3} is a natural extension of monotonic property of single variable to multivariable functions. In addition, a nonlinear function satisfying \eqref{ineq3} may not satisfy the globally Lipschitz condition, for example take $g(x)= x^3$.
  Notably, the monotone nonlinearities arise in the modeling of many physical systems like stiffening springs in constrained mechanical systems and capacitors in circuit systems \cite{fan2003observer,campbell1982singular}.
  
\noindent The first order matrix polynomial $\mu E-A$, in the indeterminate $\mu$, is called matrix pencil for the system \eqref{sys}. Moreover, system \eqref{sys} is called regular if $m=n$ and det($\mu E-A $) $\in \mathbb{R}[\mu]\setminus\{0\}$. Although, regularity ensures the existence and uniqueness of the solution to any square descriptor systems, in this paper we do not assume that \eqref{sys} is necessarily regular. Instead, we assume to have control $u(t)$ and initial state $Ex(0_{-})$ in such a way that there exists at least one solution $(x, u, v, y, z) :\mathcal{I} \rightarrow {\mathbb{R}}^{n+k+q+r+p}$ satisfying \eqref{sys}, where $\mathcal{I}=\text{dom}(x)$ is an open interval and $x \in \mathscr{C}^{1}\left( \mathcal{I} \rightarrow \mathbb{R}^n\right)$. 


		\begin{definition} \cite{berger2018observers,berger2020observers}
			Let $ \mathcal{I} \subseteq \mathbb{R}$ be an open interval. A trajectory $(x, u, v, y, z) \in \mathscr{C}( \mathcal{I} ; 
   {\mathbb{R}}^{n+k+q+r+p})$ is called solution of \eqref{sys}, if $x \in \mathscr{C}^1( \mathcal{I} ; {\mathbb{R}}^{n})$ and \eqref{sys} holds for all $t \in  \mathcal{I}$. The set 
			$$\mathscr{B}:=\{(x, u, v, y, z) \in \mathscr{C}( \mathcal{I} ; {\mathbb{R}}^{n+k+q+r+p}) \mid  (x, u, v, y, z) \; \text{is a solution of \eqref{sys}} \}$$
			of all possible solution trajectories is called the behavior of system \eqref{sys}.
		\end{definition}
		
		In order to estimate the functional state $z(t)$ in \eqref{sys}, consider the following functional filter of the form
		\begin{subequations}\label{obs}
			\begin{eqnarray}
				\dot{w}(t)&= & {N}w(t)+TBu(t)+{L}y(t)+{T}Fg(H\hat{z},u),\label{obsa}\\
				\hat{z}(t)&= & w(t)+{M}y(t),\label{obsb}
			\end{eqnarray}
		\end{subequations}
		 where $w(t) \in \mathbb{R}^{p}$ is the state vector of the filter and $\hat{z}(t) \in \mathbb{R}^{p}$ is the estimation of functional $z(t)$ in \eqref{sys}. The above matrices $N \in \mathbb{R}^{p \times p}$, $T \in \mathbb{R}^{p \times m}$, $L \in \mathbb{R}^{p \times r}$ and $M  \in \mathbb{R}^{p \times r}$  are unknown matrices.
		
		\begin{definition}\label{def2}
			System \eqref{obs} is said to be an  functional $\mathcal{H}_{\infty}$ filter for \eqref{sys}, if for every $(x, u, v, y, z) \in \mathscr{B}$, there exists $w \in \mathscr{C}^{1}\left( \mathcal{I} ; \mathbb{R}^p\right)$ and $\hat{z} \in \mathscr{C}\left( \mathcal{I} ; \mathbb{R}^p\right)$ such that $(w, u, y, \hat{z})$ satisfy \eqref{obs} for all $t \in \mathcal{I}$, and for all such $w, \hat{z}$ the following properties hold:
			\begin{enumerate}
				\item If $v \equiv 0$ and $e(t)=z(t)-\hat{z}(t)$, then $e(t) \rightarrow 0$ for $t \rightarrow \infty$.
				\item If $v$ is not identically zero and $e(t)=z(t)-\hat{z}(t)$, then  
	\begin{eqnarray*}
    \int_0^{t_f} e^T(\tau) e(\tau) \mathrm{d} \tau \leq \gamma^2 \Big(\beta + \int_0^{t_f} v^T(\tau) v(\tau) \mathrm{d} \tau\Big),  \quad t_f>0,
 \end{eqnarray*}
	where $\gamma >0$ and $\beta \geq 0$ are two real numbers.
			\end{enumerate}
		\end{definition}

In this paper, we discuss the problem of designing matrices $N$, $T$, $L$ and $M$ such that \eqref{obs} becomes a functional $\mathcal{H}_{\infty}$ filter for a given system of the form \eqref{sys}. 

We conclude this section by recalling the following fundamental result required for our analysis in the next section.

\begin{lem} \cite{piziak2007matrix}
			The linear system $\mathcal{X}\mathcal{Y} = \mathcal{Z}$ has a solution for $\mathcal{X}$ if and only if 
			$$
			\operatorname{rank}\begin{bmatrix}
				\mathcal{Y}\\
				\mathcal{Z}
			\end{bmatrix} = \operatorname{rank}\mathcal{Y}.$$
			Moreover, $$\mathcal{X}=\mathcal{Z}\mathcal{Y}^+ - \mathcal{V}(\mathcal{I}-\mathcal{Y}\mathcal{Y}^+),$$ where $\mathcal{Y}^+$ is the MP-inverse of $\mathcal{Y}$ and $\mathcal{V}$ is any matrix of appropriate dimension.
		\end{lem}

\section{Functional $\mathcal{H}_{\infty}$ filter design}\label{main}
		
\noindent We begin by formulating some linear matrix equations and solving them to determine the coefficient matrices of the functional $\mathcal{H}_{\infty}$ filter \eqref{obs}. Subsequently, we employ an LMI approach to demonstrate that the designed filter satisfies both the properties in Definition \ref{def2}. Throughout the section, we assume that \eqref{sys} satisfies a rank condition:
   \begin{eqnarray}{\label{rank}}
					\operatorname{rank}\begin{bmatrix}
					    E  & A  \\
						C  & 0  \\
						0  & C  \\
						0  & K  \\
						K  & 0 
					\end{bmatrix}
					 =
					\operatorname{rank}\begin{bmatrix}
					    E  & A \\
						C  & 0 \\
						0  & C \\
						0  & K 
					\end{bmatrix}.
				\end{eqnarray}
\textbf{Error dynamics and design of coefficient matrices for functional $\mathcal{H}_{\infty}$ filter \eqref{obs} :}
\noindent Let $e(t) = z(t)-\hat{z}(t)$ be the error between the actual and the estimated functional vector. Define $e_1(t) = TEx(t)- w(t)$. Then from \eqref{sys} and \eqref{obs}, we obtain
    \begin{eqnarray}{\label{eet}}
            e(t) &=&Kx(t)-(w(t)+My(t)), \notag \\
            &=&e_{1}(t)-(TE+MC-K)x(t)-MGv(t),
        \end{eqnarray}   
		
	\noindent and
	\begin{eqnarray}{\label{e1dot}}
\dot{e}_{1}(t) &=& {T}{E}\dot{x}(t)-\dot{w}(t), \notag \\
			&=& Ne_1(t)+(TD-LG)v(t) +\left(TA-LC-NTE\right)x(t)
	+ {T}F\big(g(Hz,u)-g(H\hat{z},u)\big), \notag \\
 & = &  Ne_1(t)+(TD-LG)v(t) +\left(TA-LC-NTE\right)x(t)+ {T}F\Delta g,
		\end{eqnarray}
   where $\Delta g  =  g(Hz,u)-g(H \hat{z},u)$. Thus, if the filter parameter matrices $N,~T,~L,$ and $M$ in \eqref{obs} satisfy the matrix equations
        	\begin{subequations}\label{cond}
					\begin{eqnarray}
						TA-LC-NTE &=& 0, {\label{conda}}\\
						TE + MC - K &=& 0, {\label{condb}}
					\end{eqnarray}
				\end{subequations}
  \noindent then, \eqref{eet} and \eqref{e1dot} infer that the error dynamics is determined by the equations
        \begin{subequations} \label{et}
		\begin{eqnarray}
         \dot{e}_1(t) &=& Ne_1(t)+(TD-LG)v(t) + {T}F\Delta g, \label{eta} \\
            e(t) &=& e_{1}(t)-MGv(t). \label{etb}
        \end{eqnarray}   
         \end{subequations}
	\noindent Notably, Eq. \eqref{conda} is nonlinear in the unknown matrices. Therefore, by substituting $TE= K-MC$ from \eqref{condb} into \eqref{conda}, and taking 
 \begin{eqnarray}\label{L}
 P  = NM-L,
 \end{eqnarray}		
 we obtain a linear matrix equation 
			\begin{eqnarray}\label{ta}
				TA+PC-NK =  0.
			\end{eqnarray}		
	Thus, we can rewrite Eq. \eqref{condb} and \eqref{ta} as
 
			\begin{eqnarray}\label{psitheta}
				\begin{bmatrix}
					T & M & P & N 
				\end{bmatrix}\Psi  =  \Theta,
			\end{eqnarray}
    where $\Psi=\begin{bmatrix}
					E  & A  \\
					C  & 0  \\
					0  & C  \\
					0  & -K 
				\end{bmatrix}$ and $\Theta= \begin{bmatrix}
					K & 0 
				\end{bmatrix}.$
    
\noindent Now, Lemma $1$ reveals that Eq. \eqref{psitheta} is solvable for the unknowns $T,~M,~P,~N$ if and only if \eqref{rank} holds, and general solution to \eqref{psitheta} is
			\begin{eqnarray}{\label{gensol}}
				\begin{bmatrix}
					T & M & P & N 
				\end{bmatrix}  =  \Theta\Psi^+-Z(I-\Psi\Psi^+),
			\end{eqnarray}
			where ${Z}$ is any matrix of suitable dimension. Therefore, we obtain
	
        \begin{subequations}\label{pmat}
     \begin{eqnarray}
				{T} &=& {T_1}-{Z}{T_2},\label{pmat1}\\
				{M} &=& {M_1}-{Z}{M_2},\label{pmat2}\\
				{P} &=& {P_1}-{Z}{P_2},\label{pmat3}\\
				{N} &=& {N_1}-{Z}{N_2},\label{pmat4}
			\end{eqnarray}
         \end{subequations}
			where
			 \begin{subequations}\label{tmpn}
     \begin{eqnarray}
				T_1 = \Theta\Psi^+\begin{bmatrix}
					I & 0 & 0 & 0
				\end{bmatrix}^T,~~T_2=(I-\Psi\Psi^+)\begin{bmatrix}
					I & 0 & 0 & 0
				\end{bmatrix}^T, \label{tmpna}\\
				M_1 = \Theta\Psi^+\begin{bmatrix}
					0 & I & 0 & 0
				\end{bmatrix}^T,~~M_2=(I-\Psi\Psi^+)\begin{bmatrix}
					0 & I & 0 & 0
				\end{bmatrix}^T, \label{tmpnb} \\
				P_1 = \Theta\Psi^+\begin{bmatrix}
					0 & 0 & I & 0
				\end{bmatrix}^T,~~P_2=(I-\Psi\Psi^+)\begin{bmatrix}
					0 & 0 & I & 0
				\end{bmatrix}^T, \label{tmpnc} \\
				N_1 = \Theta\Psi^+\begin{bmatrix}
					0 & 0 & 0 & I
				\end{bmatrix}^T, ~~N_2=(I-\Psi\Psi^+)\begin{bmatrix}
					0 & 0 & 0 & I
				\end{bmatrix}^T \label{tmpnd}.
			\end{eqnarray}
         \end{subequations}
    Thus, our remaining task is to find $Z$ in such a way that
the system \eqref{obs} with the above coefficient matrices satisfies all
conditions for functional $\mathcal{H}_{\infty}$ filter as in Definition $2$. 

\noindent Before stating the main theorem, we assume that $Z_1$ is an arbitrary matrix of the same dimension as $Z$ and set the following notations.
\begin{subequations}\label{bar}
    \begin{eqnarray}
    \bar{M_2}&=&M_2G,\label{bara}\\
    Z&=&Z_1(I-\bar{M_2}\bar{M_2}^+), \label{barb} \\
    \mathcal{T}_2&=& (I - \bar{M_2}\bar{M_2}^+)T_2, \label{barc}
    \\
    \mathcal{M}_2&=& (I - \bar{M_2}\bar{M_2}^+)M_2,\label{bard} \\
    \mathcal{P}_2&=& (I - \bar{M_2}\bar{M_2}^+)P_2, \label{bare} \\
			\mathcal{N}_2 &=& (I - \bar{M_2}\bar{M_2}^+)N_2,\label{barf} \\
   B_1&=&T_1D-N_1M_1G+P_1G, \label{barg} \\
   B_2&=&\mathcal{T}_2D-\mathcal{N}_2M_1G+\mathcal{P}_2G, \label{barh}\\
    \tilde{H} &=&M_1G, {\label{bari}} \\
    \mathscr{H} &=&H^TH. {\label{barj}}
\end{eqnarray}
\end{subequations}

 \noindent Therefore, due to the property of MP-inverse that $\bar{M_2}\bar{M_2}^+\bar{M_2}=\bar{M_2}$, we obtain
 \begin{eqnarray}\label{mp}
      \mathcal{M}_2G&=& (I - \bar{M_2}\bar{M_2}^+)\bar{M_2}=0.
 \end{eqnarray}
We are now ready to state the main theorem.
			
   \begin{theorem}
   Consider a system \eqref{sys} which satisfies the rank condition \eqref{rank} and the nonlinear function $g$ satisfies \eqref{ineq2}. Then, for a given $\gamma>0 $,   \eqref{obs} is a functional $\mathcal{H}_{\infty}$ filter with system coefficient matrices \eqref{pmat} and error dynamics \eqref{et} if there exist two matrices $ Q = Q^T > 0$, and $Y$ such that the following LMI holds.
				\begin{eqnarray}{\label{lmi}}
     \tilde{\Pi} = \begin{bmatrix}
			\mathcal{A}_{11}+I & \mathcal{A}_{12}  &  \mathcal{A}_{13}-\tilde{H}\\
				\mathcal{A}_{12}^{T} & 0 &  \mathcal{A}_{23}\\
						(\mathcal{A}_{13}-\tilde{H})^{T} & \mathcal{A}_{23}^{T} &   \mathcal{A}_{33}+\tilde{H}^T\tilde{H}-\gamma^2 I
					\end{bmatrix} & \leq & 0,
				\end{eqnarray}
				where 
    \begin{subequations}\label{cof}
        \begin{eqnarray}
					\mathcal{A}_{11} &=& N_{1}^T Q+Q N_{1}-\mathcal{N}_{2}^{T} Y^T-Y\mathcal{N}_{2} -\rho \mathscr{H}, {\label{cof1}}\\
					\mathcal{A}_{12} &=& (QT_{1}-Y\mathcal{T}_{2})F+H^T,{\label{cof2}}\\
                     \mathcal{A}_{13} &=& Q B_{1}-Y B_{2} + \rho \mathscr{H}\tilde{H},{\label{cof3}}\\
                     \mathcal{A}_{23} &=& -H\tilde{H},{\label{cof4}} \\
					\mathcal{A}_{33} &=& -\rho\tilde{H}^{T}\mathscr{H}\tilde{H},{\label{cof5}}\\
                     Y &=& QZ_1.{\label{cof6}}
				\end{eqnarray}
            \end{subequations}
				
			\end{theorem}
			
           \begin{proof}
           We break the proof into three steps.
           
           \noindent \textbf{Step 1:}
          In this step, we transform the error dynamics \eqref{et} in a convenient form so that solution $(Q,Y)$ of LMI \eqref{lmi} provides $Z_1= Q^{-1}Y$ (cf. \eqref{cof6}) in such a way that \eqref{obs} becomes a functional $\mathcal{H}_{\infty}$ filter of \eqref{sys} with the coefficient matrices \eqref{pmat}. Clearly, \eqref{bara}-\eqref{bari} infer that \eqref{pmat} can be rewrite as
      \begin{subequations} \label{pm}
     \begin{eqnarray}
			    T &=& T_1-Z_1\mathcal{T}_2,{\label{pm1}}\\
				M &=& M_1-Z_1\mathcal{M}_2,{\label{pm2}}\\
				P &=& P_1-Z_1\mathcal{P}_2,{\label{pm3}}\\
				N &= &N_1-Z_1\mathcal{N}_2.{\label{pm4}} 
			\end{eqnarray}
    \end{subequations}
   Moreover, if we denote
   \begin{eqnarray}
       \mathcal{B} &= & B_1-Z_1B_2, {\label{pm5}}
   \end{eqnarray}
    
			\noindent then, by using \eqref{mp} and \eqref{pm5}, we obtain
   \begin{subequations}
       \begin{eqnarray}\label{td}
				TD-LG &=& TD - (NM-P)G  \notag\\
    &=& (T_1-Z_1\mathcal{T}_2)D - [(N_1-Z_1\mathcal{N}_2)(M_1-Z_1\mathcal{M}_2)-(P_1-Z_1\mathcal{P}_2)]G \notag \\
    &=& B_1-Z_1B_2 = \mathcal{B}, \label{tda}\\
    \text{and} \qquad
    MG &=& (M_1-Z_1\mathcal{M}_2)G=M_1G \label{tdb}.
			\end{eqnarray}
   \end{subequations}
   
  \noindent Thus, we can rewrite the error dynamics \eqref{et} as
	\begin{subequations}\label{newerr}
		\begin{eqnarray}
			\dot{e}_{1}(t) &=& Ne_{1}(t) + TF\Delta g+ \mathcal{B}v(t) ,\label{newerra}\\
			e(t) &=& {e}_{1}(t)-\tilde{H}v(t), \label{newerrb}
				\end{eqnarray}
			\end{subequations}
   where $N$, $T$, $\mathcal{B}$ and $\tilde{H}$ are defined in \eqref{pm4}, \eqref{pm1}, \eqref{pm5} and \eqref{bari} respectively.
			
  \noindent Since $g$ satisfies the generalized monotone property, using \eqref{newerrb}, Eq. \eqref{ineq2} reduces to 
    \begin{eqnarray} \label{secprop}
       \Delta x^{T}\Delta g + \Delta g^{T}\Delta x & \geq& \rho \|\Delta x\|^2 ,\notag \\
       &=& \rho  \|He(t)\|^2 ,\notag \\
        &=& \rho  \|H\big({e}_{1}(t)-\tilde{H}v(t)\big)\|^2 ,\notag \\
		&=& \rho \big(e_1^{T}(t)\mathscr{H}e_1(t)- e_1^{T}(t)\mathscr{H}\tilde{H}v(t) 
		- v^{T}(t)\tilde{H}^{T}\mathscr{H}e_1(t) \notag\\
          & & +v^{T}(t)\tilde{H}^{T}\mathscr{H}\tilde{H}v(t)\big).
    \end{eqnarray}
   \noindent \textbf{Step 2:} In this step, we show that if LMI \eqref{lmi} holds, then the error dynamics \eqref{newerr} satisfies condition $(b)$ in Definition \ref{def2}. Consider the Lyapunov function $V$ such that
   \begin{eqnarray*}
       V : \mathbb{R}^p \rightarrow \mathbb{R}, \; e_1(t) \rightarrow e_{1}^T(t)Qe_{1}(t),
   \end{eqnarray*}
   where $Q=Q^{T}>0$ and $V(e_1(0)) = \gamma^2 \beta$ for some $\beta \geq 0$.
   Thus, by using \eqref{newerr}, we obtain
			\begin{eqnarray*}\label{leqq} 
				    \dot{V}(e_1(t)) &=& \dot{e_{1}}^{T}(t) Q e_{1}(t)+e_{1}^T(t) Q \dot{e_{1}}(t) \notag \\
					&=&(N e_{1}(t)+T F\Delta g+\mathcal{B}v(t))^{T} Q e_{1}(t)+e_{1}^{T}(t) Q(N e_{1}(t)+TF \Delta g+\mathcal{B}v(t)) \notag \\
					&= &e_{1}^{T}(t)\big((N_{1}-Z_1\mathcal{N}_{2})^{T}Q + Q(N_{1}-Z_1\mathcal{N}_{2})\big)e_{1}(t) +e_{1}^{T}(t)\big(Q (T_{1}-Z_1\mathcal{T}_{2})\big)F \Delta g \notag \\
					&& +\Delta g^{T}\big(Q \left(T_{1}-Z_1\mathcal{T}_{2}\right)F\big)^{T} e_{1}(t)
					+e_{1}^T(t) Q (B_1-Z_1 B_2) v(t) \\
                    &&+v^T(t) \big(Q(B_1-Z_1B_2)\big)^T e_{1}(t) 
     \end{eqnarray*}
   Therefore, \eqref{secprop} reveals that
   \begin{eqnarray}\label{leq}
                 \dot{V}(e_1(t)) &{\leq}& {e_{1}}^{T}(t)\big((N_{1}-Z_1\mathcal{N}_{2})^{T} Q+Q(N_{1}-Z_1\mathcal{N}_{2})\big)e_1(t)+e_{1}^{T}(t)\big(Q \left(T_{1}-Z_1\mathcal{T}_{2}\right)\big)F \Delta g \notag \\
				&& +\Delta g^{T}\big(Q(T_{1}-Z_1\mathcal{T}_{2})F\big)^{T}e_{1}(t) + e_{1}^T(t) Q (B_1-Z_1 B_2) v(t) \notag \\
                &&+v^T(t) \big(Q(B_1-Z_1B_2)\big)^T e_{1}(t) +\big(H({e}_{1}(t) -\tilde{H}v(t))\big)^{T} \Delta g \notag \\
			    && + \Delta g^{T}\big(H({e}_{1}(t)-\tilde{H}v(t))\big) -\rho e_1^{T}(t)\mathscr{H}e_1(t)+ \rho e_1^{T}(t)\mathscr{H}\tilde{H}v(t) \notag \\
					&& + \rho v^{T}(t)\tilde{H}^{T}\mathscr{H}e_1(t) - \rho v^{T}(t)\tilde{H}^{T}\mathscr{H}\tilde{H}v(t) \notag \\
					&=&{e_{1}}^{T}(t)\big((N_{1}-Z_1\mathcal{N}_{2})^{T} Q+Q(N_{1}-Z_1\mathcal{N}_{2})-\rho \mathscr{H}\big)e_1(t) \notag \\
					&& +e_{1}^{T}(t)\big(Q (T_{1}-Z_1\mathcal{T}_{2})F+H^{T}\big)\Delta g +\Delta g^{T}\big(\big(Q (T_{1}-Z_1\mathcal{T}_{2})F\big)^{T}+H\big)e_{1}(t) \notag \\
				&&  +e_{1}^T(t) \big(Q (B_1-Z_1 B_2)+ \rho \mathscr{H}\tilde{H}\big)v(t) +v^T(t) \big(\big(Q(B_1-Z_1B_2)\big)^T +\rho \tilde{H}^{T}\mathscr{H}\big) e_{1}(t) \notag \\
					&& - v^{T}(t)\rho\tilde{H}^{T}\mathscr{H}\tilde{H}v(t)-v^{T}(t)\tilde{H}^{T}H^T\Delta g - \Delta g^{T} H\tilde{H}v(t),
				\end{eqnarray}
				which is equivalent to the fact that
				\begin{eqnarray}{\label{vdot}}
                    \dot{V}(e_1(t))&\leq&	\zeta^{T}(t)\begin{bmatrix}
						\mathcal{A}_{11} & \mathcal{A}_{12}  &  \mathcal{A}_{13}\\
						\mathcal{A}_{12}^{T} & 0 &  \mathcal{A}_{23}\\
						\mathcal{A}_{13}^{T} & \mathcal{A}_{23}^{T} &  \mathcal{A}_{33}
					\end{bmatrix}\zeta(t),\quad
   \text{where}\; \zeta(t) =\begin{bmatrix}
					e_1(t) \\
					\Delta g\\
                    v(t)
				\end{bmatrix}.
				\end{eqnarray}
Now, from \eqref{newerrb} and \eqref{vdot}, we obtain the following inequality
    \begin{eqnarray*}
        \dot{V}(e_1(t)) +e^T(t)e(t)-\gamma^2v^T(t)v(t)&\leq& \zeta^T(t)\tilde{\Pi}\zeta(t),
    \end{eqnarray*}
  where $\tilde{\Pi}$ is the same as in \eqref{lmi}.\\
Thus, if $\tilde{\Pi} \leq 0$, then the above inequality gives that
\begin{eqnarray*}
   \dot{V}(e_1(t)) &\leq& \gamma^2 v^T(t) v(t)-e^T(t) e(t).
\end{eqnarray*}
    
\noindent Thus, by integrating both sides of the above inequality from $0$ to $t_f$, we obtain
				\begin{eqnarray*}
				    \int_0^{t_f} \dot{V}(\tau) \mathrm{d} \tau & \leq & \int_0^{t_f} \gamma^2 v^T(\tau) v(\tau) \mathrm{d} \tau-\int_0^{t_f} e^T(\tau) e(\tau) \mathrm{d} \tau,
				\end{eqnarray*}
				or equivalently 
    \begin{eqnarray*}
         V(e_1(t_f))-V(e_1(0)) \leq \int_0^{t_f} \gamma^2 v^T(\tau) v(\tau) \mathrm{d} \tau-\int_0^{t_f} e^T(\tau) e(\tau) \mathrm{d} \tau.
    \end{eqnarray*}
 Since $V(e_1(t_f)) \geq 0$, we conclude that 
 \begin{eqnarray*}
    \int_0^{t_f} e^T(\tau) e(\tau) \mathrm{d} \tau \leq \gamma^2 \Big( \beta +  \int_0^{t_f} v^T(\tau) v(\tau) \mathrm{d} \tau \Big),
 \end{eqnarray*}
  which infers that the condition $2.$ in Definition $2$ is satisfied if \eqref{lmi} holds.
 
\noindent \textbf{Step 3:} This step proves that, if $v \equiv 0$, the error dynamics \eqref{newerr} also satisfies condition $1.$ in Definition \ref{def2}. Take $v \equiv 0$, then the error dynamics \eqref{newerr} reveals that $e(t)=e_1(t)$ and 
		\begin{eqnarray}\label{newerr1}
			\dot{e_1}(t) &=& ({N_{1}-Z_1\mathcal{N}_{2})}e_{1}(t) + ({T_{1}-Z_1\mathcal{T}_{2}})F\Delta g,
				\end{eqnarray}
      and Eq. \eqref{leq} reduces to
     \begin{eqnarray*}\label{leq1} 
			 \dot{V}  & \leq & {e_{1}}^{T}(t)\big((N_{1}-Z_1\mathcal{N}_{2})^{T} Q+Q(N_{1}-Z_1\mathcal{N}_{2})-\rho \mathscr{H}\big)e_1(t)+e_{1}^{T}(t)\big(Q \left(T_{1}-Z_1\mathcal{T}_{2}\right)F \notag \\
			&& +H^T\big)\Delta g +\Delta g^{T}\big((Q \left(T_{1}-Z_1\mathcal{T}_{2}\right)F)^{T}+H\big)e_{1}(t), \notag \\
                        & = &  \begin{bmatrix}
                            e_1(t) \\
                            \Delta g
                        \end{bmatrix}^T
                       \begin{bmatrix}
					    N_{1}^T Q+Q N_{1}-\mathcal{N}_{2}^{T} Y^T-Y\mathcal{N}_{2} -\rho \mathscr{H}  &  (QT_{1}-Y\mathcal{T}_{2})F+H^T\\
						\big((QT_{1}-Y\mathcal{T}_{2})F\big)^{T}+H  &  0
					\end{bmatrix}  \begin{bmatrix}
                            e_1(t) \\
                            \Delta g
                        \end{bmatrix}.
                 \end{eqnarray*}
 Thus, $\dot{V}<0$ if
				\begin{equation}\label{lmi1}
					\Omega = \begin{bmatrix}
					    N_{1}^T Q+Q N_{1}-\mathcal{N}_{2}^{T} Y^T-Y\mathcal{N}_{2} -\rho \mathscr{H}+I  &  (QT_{1}-Y\mathcal{T}_{2})F+H^T\\
						\big((QT_{1}-Y\mathcal{T}_{2})F\big)^{T}+H  &  0
					\end{bmatrix} \leq 0
				\end{equation}
              which is a principle submatrix of $\tilde{\Pi}$. Thus, if \eqref{lmi} holds, the error dynamics \eqref{newerr1}  is asymptotically stable \cite[Thm. $5.14$]{barnett1985introduction} and hence $e(t)=e_1(t) \rightarrow 0$ as $t \rightarrow \infty$.
              \end{proof}
              
            
			\noindent Based on Theorem $1$, we now summarize the functional $\mathcal{H}_{\infty}$ filter design procedure in Algorithm $1$ below.
			\FloatBarrier
			\begin{algorithm}\caption{Computational steps to construct functional $\mathcal{H}_{\infty}$ filter \eqref{obs} for system \eqref{sys}. }
				\textbf{Step $1$:} Check whether \eqref{rank} is satisfied. If yes, go to the next step.\\
				\textbf{Step $2$:}  Find $T_i$, $M_i$, $P_i$ and $N_i$ for $i=1,2$ from \eqref{tmpna} to \eqref{tmpnd}. \\
				\textbf{Step $3$:} Determine $\mathcal{T}_2$, $\mathcal{M}_2$, $\mathcal{P}_2$ and $\mathcal{N}_2$ from \eqref{barc} to \eqref{barf}. \\
				\textbf{Step $4$:} Compute $B_1$, $B_2$ and $\tilde{H}$ from \eqref{barg}, \eqref{barh} and \eqref{bari}, respectively. \\
				\textbf{Step $5$:}  Solve \eqref{lmi} for $Q$ and $Y$. \\
				\textbf{Step $6$:}  Compute $Z_1 = Q^{-1}Y$ and $Z = Z_1(I-\bar{M_2}\bar{M_2}^{+})$. \\
				\textbf{Step $7$:} Calculate the matrices $T$, $M$, $P$ and $N$ from \eqref{pm1} to \eqref{pm4}.\\
				\textbf{Step $8$:} Compute $L$ from \eqref{L}.
			\end{algorithm}
			\FloatBarrier	
		
\section{Illustrative Example}\label{num}
			\noindent In this section, we implement the proposed functional $\mathcal{H}_{\infty}$ filter design algorithm on a disc of mass $m$ rolling on a surface without slipping. The disc is connected to a fixed wall by a nonlinear spring with coefficients $k_1>0$, $k_2>0$ and a linear damper with damping coefficient $b>0$, as shown in Fig.$1$ \cite[pp. 78]{sjoberg2006some}. The position of the center of the disc along the surface is given by $x_1$, while $x_2$ is the translational velocity of the same point. Moreover, the angular velocity of the disc is denoted by $x_3$. A control input torque $u$ is applied at the center of the disc to rotate it, and the frictional force $\lambda$ that resists the motion can be considered an unknown bounded disturbance.   
				\begin{figure}
					\centering
					\includegraphics[height = 1.6in]{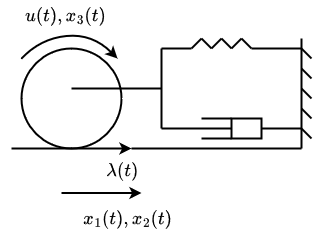}
					\caption{Rolling disc}
				\end{figure}

\noindent The system dynamics along with the kinematic constraints on motion are determined as \cite{sjoberg2006some}
\allowdisplaybreaks
\begin{equation}{\label{exp2}}
\begin{aligned}
\dot{x}_1 &=x_2, \\
						\dot{x}_2 &=-\frac{k_1}{m} x_1-\frac{k_2}{m} x_1^3-\frac{b}{m} x_2+\frac{\lambda}{m} ,\\
						0 &=x_2-r x_3 ,\\
						0 &=-\frac{k_2}{m} x_1^3-\frac{k_1}{m} x_1-\frac{b}{m} x_2+\left(\frac{1}{m}+\frac{r^2}{J}\right) \lambda-\frac{r}{J} u,
					\end{aligned}
				\end{equation}
                 where the moment of inertia about the center of the disc is $J$ and $r$ is the radius of the disc.
                 
                 If the aim is to estimate the position of the disc, then $z=x_1$, and the system \eqref{exp2} can be expressed in the form of \eqref{sys} with the following coefficient matrices:
                  \\~\\
					\noindent $E=\begin{bmatrix}
						1 & 0 & 0\\
						0 & 1 & 0 \\
						0 & 0 & 0 \\
					\end{bmatrix}$, $A=\begin{bmatrix}
						0 & 1 & 0\\
						-{k_1}/{m} & -{b}/{m} & 0 \\
						-{k_1}/{m} & 0 & -{rb}/{m} \\
					\end{bmatrix}$,
					$B=\begin{bmatrix}
						0 \\
						0 \\
						-r/J
					\end{bmatrix}$, 
					$D=\begin{bmatrix}
						0 \\
						1/m \\ 
						1/m + r^2/J
					\end{bmatrix}$, \\~\\  $F=\begin{bmatrix}
						0 \\
						-{k_2}/m \\ 
						-{k_2}/m
					\end{bmatrix}$, 
					$C=\begin{bmatrix}
						1 & 1  & 0\\
						0 & 0 & 1
					\end{bmatrix}$, $K=\begin{bmatrix}
						1 & 0 & 0 \\
					\end{bmatrix}$, $G=\begin{bmatrix}
						0.35 \\
						0.11
					\end{bmatrix}$, $H=1$.
     
\noindent     Clearly, the nonlinear function  $g = x_1^3$ satisfies \eqref{ineq} with $\rho=0$.
				For simulation purpose, we take
				$\dfrac{k_1}{m}=\dfrac{k_2}{m}=1, \; \dfrac{b}{m}=2, \; r=2, \; m=1 \; \text{and} \; J=4.$ 
			   It is now easy to check that the system satisfies the rank condition \eqref{rank}. Now, for $\gamma=1.4$ and $\beta=0.1$, we solve the LMI \eqref{lmi} by using \emph{feasp} function in MATLAB and obtain the parameter matrix $Z$ as
            \begin{center}
                $Z=\begin{bmatrix}
                    442.8408 & 136.8489 & -11.2463 & 546.9125 &  11.7694   & 38.4020 &   7.0763 &  19.2352

                \end{bmatrix}$
            \end{center}
     Using the above matrix $Z$, we compute the filter coefficient matrices by using Algorithm $1$ as
			\begin{center}
				$T=\begin{bmatrix}
					 1.0000  & -0.0000  &  0.0653 \\
				\end{bmatrix}$, 
				$L=\begin{bmatrix}
					1.0000  & -0.2614 \\
				\end{bmatrix}$, 
				$M=\begin{bmatrix}
					-0.0779  &  0.1553 \\
				\end{bmatrix} \times 1.0e{-13}$,\\
				$N=-1.0653$ and $P=\begin{bmatrix}
					 -1.0000  &  0.2614 \\
				\end{bmatrix}$.
			\end{center}
   
We plotted the true and estimated values of $z$ in Fig.$3$ by taking $x(0) = [0.1 \hspace{0.5em}  0.2 \hspace{0.5em}0.15]^{T}$, $w(0) = 0.3$, $u(t)= 0.2\sin(\pi t)$ and the random disturbance vector $\lambda(t)$, which is shown in Fig.$2$, 
     is generated by $\text{rand}$ function in MATLAB. 
            
             \begin{figure}
             \centering
	\includegraphics[width=0.52\linewidth]{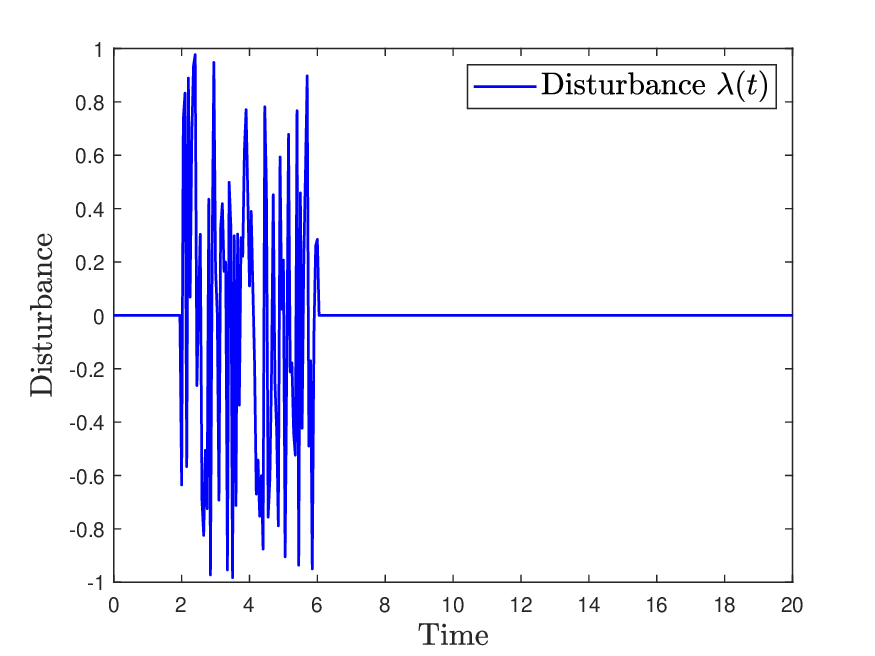}
	\caption{Random disturbance $\lambda(t) = 2 \times \text{rand}(\text{size}(t))-1;\, t \in [2s, 6s]$}
        \end{figure}
    
    \begin{figure}[H]
	\begin{subfigure}[h]{0.52\textwidth}
	\includegraphics[width=\linewidth]{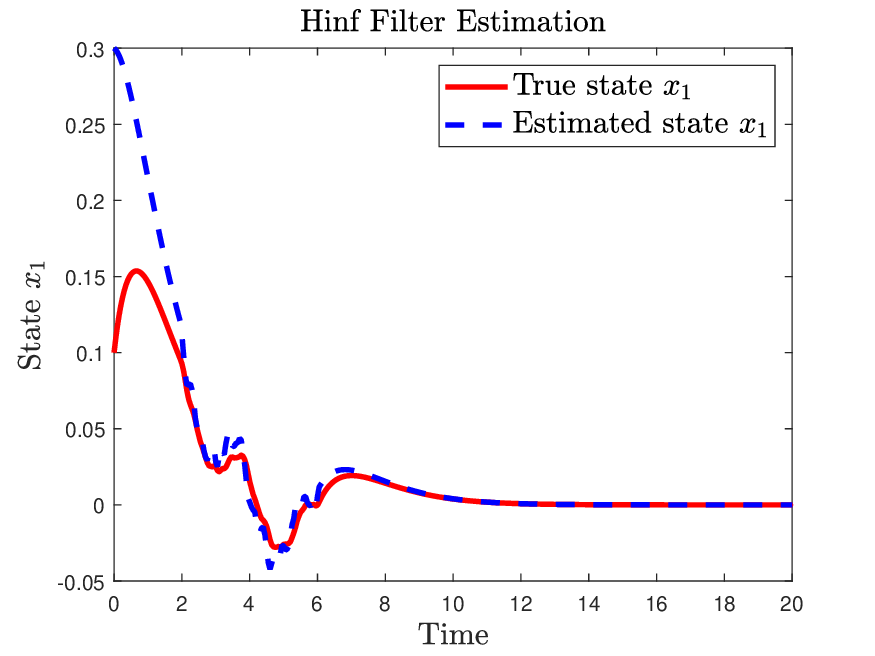}
	\caption{Time responses of true and estimated state $z = x_1$}
        \end{subfigure}
	\hfill
	\begin{subfigure}[H]{0.52\textwidth}
	\includegraphics[width=\linewidth]{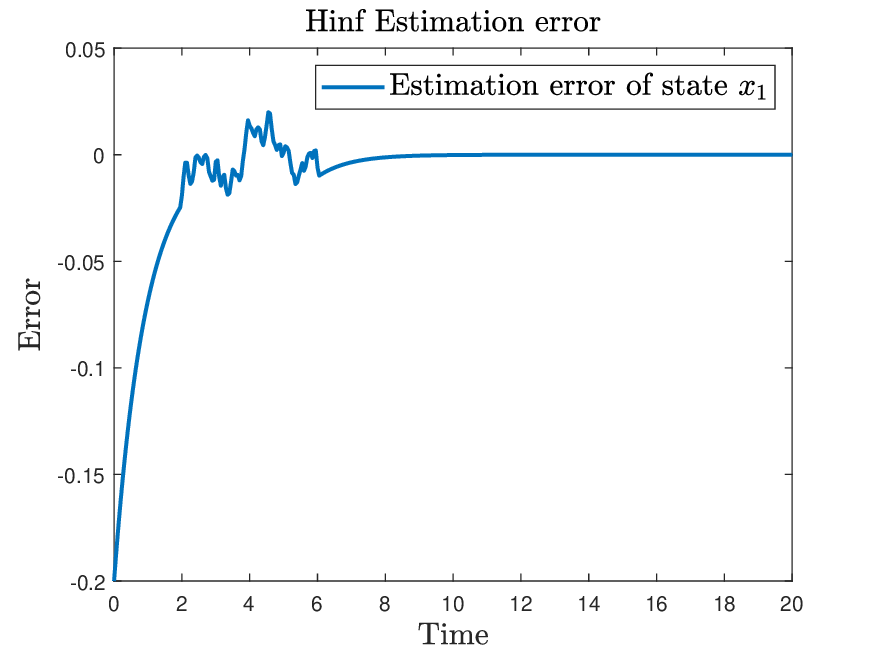}
	\caption{Time Response of error $e(t)=z(t)-\hat{z}(t)$}
	\end{subfigure}
	\caption{Functional $\mathcal{H}_{\infty}$ filter performance with bounded disturbance $\lambda(t)$}
	\end{figure}

\section{Conclusion}\label{concl}
This paper has contributed to the development of functional $\mathcal{H}_\infty$ filtering for nonlinear descriptor systems with generalized monotone nonlinearities. The existence conditions for the functional $\mathcal{H}_\infty$ filter have been established through a rank condition on system matrices and an LMI. Notably, our designed functional filter does not impose specific conditions on the matrix pencil related to the given descriptor systems, except for the essential requirement of system well-posedness. Furthermore, the designed filter exhibits robustness against unknown disturbances, as demonstrated through a numerical example.  Finally, this paper opens up future opportunities for advancing the theory and application of functional $\mathcal{H}_\infty$ filtering in nonlinear descriptor systems. For instance, future research could explore the development of conditions and algorithms for constructing the smallest possible order functional filter. Additionally, investigating the incorporation of other types of nonlinearities in the system and output presents a fascinating avenue for further exploration.

\section*{Declarations}
 \noindent \textbf{Funding information:} Rishabh Sharma acknowledges funding by the Council of Scientific and Industrial Research, New Delhi, India, for the award of SRF through grant number 09/1023(0039)/2020-EMR-I. Nutan K. Tomar acknowledges funding by the Science and Engineering Research Board, New Delhi, via Grant number CRG/2023/008861.\\
 
 \noindent \textbf{Author contributions:} Rishabh Sharma: conceptualization, writing- original draft, software, investigation, validation and writing- review and editing. Mahendra K. Gupta: Validation and writing, structuring the draft and editing. Nutan K. Tomar: conceptualization, supervision, writing-original draft, validation and review and editing.  \\
 
 \noindent \textbf{Conflict of interest:} The authors declare that they have no known competing financial interests or personal relationships that could have appeared to influence the work reported in this paper.\\
 
 \noindent \textbf{Ethical approval:} This research did not require ethical approval.\\
 
 \noindent \textbf{Data availability statement:} Data sharing is not applicable to this article as no new data were created or analyzed in this body.


\bibliography{bibfile}

\end{document}